\documentclass[10pt]{article}

\usepackage{amsmath,amssymb,epsfig,calc,url,theorem,chngcntr}
\newcommand{\proofend}{\hfill $\blacklozenge$}

\newcommand{\Z}{\mathbb{Z}}

\newcommand{\F}{\mathbb{F}}

\newcommand{\Zm}{\Z/m\Z}

\newcounter{romancnt}
\newenvironment{romanlist}
  {\begin{list}{\makebox[\labelwidth][r]{{\rm (\roman{romancnt})}}}
               {\usecounter{romancnt}%
                \topsep         0pt%
                \parskip        0pt%
                \partopsep      0pt%
                \itemsep        0pt%
                \parsep         0pt}}
  {\end{list}\vspace{-\parskip}}

  {\begin{enumerate}
                \topsep         0pt%
                \parskip        0pt%
                \partopsep      0pt%
                \parsep         0pt}
  {\end{enumerate}}

\counterwithin{equation}{section}

\newtheorem{theorem}[equation]{Theorem}

\newtheorem{lemma}[equation]{Lemma}
\newtheorem{prop}[equation]{Proposition}
\newtheorem{cor}[equation]{Corollary}

{\theorembodyfont{\rmfamily}
 \newtheorem{prob}[equation]{Problem}}
{\theorembodyfont{\rmfamily}
 }
{\theorembodyfont{\rmfamily}
 \newtheorem{defi}[equation]{Definition}}

\begin{document}

\title{Exact solutions to Waring's problem for finite fields}
\author{by\\
        {\small }\\
	Arne Winterhof\\
	and\\
        Christiaan van de Woestijne}
\date{2 October, 2008}

\maketitle

{ 
\renewcommand{\thefootnote}{}
\footnotetext[0]{1991 {\it Mathematics Subject Classification}: Primary 
                 11P05, Secondary 11T41, 90C10, 94B65.}
}

%

\section{Introduction and results}

Let $q=p^r$ be a power of a prime $p$ and denote by $\F_q$ the finite field of $q$ elements.
For a positive integer $k$, 
\emph{Waring's problem} for $\F_q$ is
the question how many summands are maximally needed to express any given
element $a$ of $\F_q$ in the form
\begin{equation} \label{EqRep}
  a = \sum_{i=1}^g x_i^k
\end{equation}
with $x_i\in\F_q$, i.e., as a sum of $k$th powers of elements of $\F_q$. We can
then define the \emph{Waring function} $g(k,q)$ as
the maximal number of summands needed to express all elements of $\F_q$ as sums of $k$th
powers.

We note that, by an easy argument, we have $g(k,q)=g(k',q)$, where $k'=
\gcd(k,q-1)$. Hence, we will assume from now on that $k$ divides $q-1$.

Several authors have established bounds on the value of $g(k,q)$ for various 
choices of the parameters $k$ and $q$ -- a survey is given in
\cite{wi1}. For the cases where the exponent $k$ is small compared to $q$,
there are strong results. For example, whenever $2\le k<q^{1/4}+1$, it follows that
$g(k,q)= 2$ by a direct application of the Weil bound for the number of
points on varieties over finite fields \cite{Small,Weil,wi1}.

In this paper, we will look at the cases where the exponent $k$ is \emph{large}
compared to $q$, and we will obtain not only a bound, but the \emph{exact
value} of $g(k,q)$ for two infinite families of pairs $(k,q)$. Our main results
are the following.

\begin{theorem}
\label{Thm1}
  Let $p$ and $r$ be primes such that $p$ is a primitive root modulo~$r$.
  Then we have
  $$
    g\left(\frac{p^{r-1}-1}{r},p^{r-1}\right)=\frac{(p-1)(r-1)}{2}.
  $$
\end{theorem}

\begin{theorem}
\label{Thm2}
  Let $p$ and $r$ be odd primes such that $p$ is a primitive root modulo~$r$.
  Then we have
  $$
    g\left(\frac{p^{r-1}-1}{2r},p^{r-1}\right)
        = \begin{cases}
             \lfloor \frac{pr}4 - \frac{p}{4r} \rfloor & \text{ if } r < p; \\
             \lfloor \frac{pr}4 - \frac{r}{4p} \rfloor & \text{ if } r \ge p.
          \end{cases}
  $$
\end{theorem}

{\em Remarks}. 1. Theorem~\ref{Thm1} improves the lower bound of
\cite[Theorem~2]{wi2}.

2. The value
$$
  g(p-1,p)=p-1
$$
can be regarded as complement of Theorem~\ref{Thm1} in the case that $r=1$.

3.
The values
\begin{align*}
  g((p-1)/2,p)&=(p-1)/2 \\
  \intertext{and}
  g((p^2-1)/4,p^2)&=p-1 \quad \mbox{ if }p\equiv 3 \bmod 4
\end{align*}
can be regarded as complements of Theorem~\ref{Thm2} in the case that
$r=1$ or $r=2$, respectively.

~

The proofs of our results rest on the resolution (Theorems \ref{ThmG} and
\ref{ThmMR}) of two instances of a combinatorial problem, which will be given
in detail in the next section. The problem may be loosely formulated as the
determination of the covering radius of cyclic codes in the so-called \emph{Lee
metric} (instead of in the usual Hamming metric). There is also a connection to
the determination of the diameter of \emph{Waring graphs} in graph theory
\cite{gaso}.

Section \ref{SecG} is devoted to the proof of Theorem \ref{ThmG}, from which
Theorem \ref{Thm1} follows.

The proof of Theorem \ref{ThmMR}, which implies Theorem \ref{Thm2}, is
much more involved. We prove that the values given in this theorem are
\emph{upper bounds} for the Waring function in Section \ref{SecUpper}, while in
Section \ref{SecCon} we show that the bounds are attained. Everything is put
together in Section \ref{SecProof}. The proof is constructive, in the sense
that it gives an algorithm to construct elements in $\F_q$ that need a maximal
number of terms to express them as sums of $k$th powers. 

An implementation of
this algorithm using the KASH computer algebra system (version 2.x) is
available from the second author's homepage \cite{CvdWImpl}.

\subsection*{Acknowledgements}

We want to thank Hendrik Lenstra for suggesting this way of attacking Waring's
problem.

The research that led to this publication was supported by the Austrian Science
Foundation FWF, in Linz by Projects S8313 and P19004-N18, and in Graz by
Project S9606, which is part of the Austrian National Research Network
``Analytic Combinatorics and Probabilistic Number Theory.''

\section{A combinatorial reformulation} \label{SecComb}

Let $m$ and $r$ be positive integers, and consider the free $\Z/m\Z$-module
$$
  V = \left( \Z / m\Z \right)^r.
$$
Let $g_1,\ldots,g_r$ be a basis of $V$, and define $V'$ as the quotient of $V$
by the relation $g_1+\ldots+g_r=0$. Then every element $v$ of $V'$ has multiple
representations
\begin{equation} \label{EqModuleRep}
  v = \sum_{i=1}^r v_ig_i \qquad (v_i\in \Z/m\Z),
\end{equation}
and one is interested in the size of the most economical representation.
Here, ``economical'' of course must be defined, and we will do this in two
distinct ways.

The first definition that we use assigns to each element $x$ of $\Z/m\Z$ its
least residue modulo $m$, denoted by $\bar{x}\in \{0,1,\ldots,m-1\}$, and looks
at
$$
  \|(v_1,\ldots,v_r)\|_1 =_{\rm def} \sum_{i=1}^r \bar{v_i}.
$$
The second uses the absolute least residue modulo $m$,
$$
  |x|=\min\{\bar{x},m-\bar{x}\},
$$
and looks at the {\em Lee norm}
$$
  \|(v_1,\ldots,v_r)\|_2 =_{\rm def} \sum_{i=1}^r |v_i|.
$$
It is clear that if the coefficients $(v_1,\ldots,v_r)$ and
$(v_1',\ldots,v_r')$ both represent the same element $v$ in the form 
\eqref{EqModuleRep}, then we have
$$
  (v_1',\ldots,v_r') = (v_1,\ldots,v_r) + x {\bf e}
$$
for some $x\in\Z/m\Z$, where ${\bf e}$ denotes the vector $(1,1,\ldots,1)$. 

We now give the precise definition of ``economic''. We call a vector in $V$
\emph{admissible} if 
$$
  \|{\bf v}\|_i \le \|{\bf v}+x\cdot {\bf e}\|_i \quad \text{ for all }
  x\in\Z/m\Z,  
$$
where $i$ is either $1$ or $2$, depending on the context. The problem to be 
solved is the following, where ``norm'' is one of $\|\cdot\|_1$ or
$\|\cdot\|_2$.

\begin{prob} \label{prob}
  Given positive integers $m$ and $r$, what is the largest possible norm
  of an admissible vector in $(\Z/m\Z)^r$?
\end{prob}

We will provide a complete answer to this question. Define the \emph{norm bound
functions} $g(m,r)$ and $h(m,r)$ for positive integers $m$ and $r$ by 
\begin{align} 
    g(m,r) &= \frac{mr - m - r +\gcd(m,r)}{2}; \\
    \label{EqDefH}
    h(m,r) &= 
    \begin{cases}
               \frac{mr}4 & \text{ if $m$ and $r$ are even}; \\
	       \lfloor \frac{mr}4 - \frac12 \rfloor & 
	          \text{ if $m$ is even, $r$ is odd, and $r > m$}; \\
	       \lfloor\frac{mr}4 - \frac{r}{4m} \rfloor &
	          \text{ if $m$ is odd and $r> m$}; \\
	       \lfloor \frac{mr}4 - \frac12\rfloor & 
	          \text{ if $m$ is odd, $r$ is even, and $r < m$}; \\
	       \lfloor \frac{mr}4 - \frac{m}{4r}\rfloor &
	          \text{ if $r$ is odd and $r \le m$}.
	     \end{cases}
\end{align}

Note that $g(m,r)$ is always an integer.

\begin{theorem} \label{ThmG}
Let $m$ and $r$ be positive integers, and let ${\bf v}$ be an admissible
vector in $V=(\Zm)^r$ of maximal norm $\|{\bf v}\|_1$. Then
$$
  \|{\bf v}\|_1 = g(m,r).
$$
\end{theorem}

\begin{theorem} \label{ThmMR}
Let $m$ and $r$ be positive integers, and let ${\bf v}$ be an admissible
vector in $V=(\Z/m\Z)^r$ of maximal Lee norm $\|{\bf v}\|_2$. Then 
$$
  \|{\bf v}\|_2 = h(m,r).
$$
\end{theorem}

See the next sections for the proofs of these results.

We note that Problem \ref{prob} given above can be reinterpreted in terms of
\emph{covering radii} of linear codes, with respect to the \emph{Lee metric}.
This link was also observed by Helleseth in \cite{he}.

The covering radius is a fundamental parameter of a code and has extensively
been studied. For example the subject is treated in the 
survey~\cite{brlipl} and in the monograph~\cite{coholilo}.
Let $C\subseteq (\Zm)^r$ be a code over $\Zm$ of length $r$. 
We say that a vector is $\rho$-covered by a code if it has {\em Lee-distance}
at most $\rho$ from at least one codeword. (The Lee distance of
$(a_1,\ldots,a_r)$, $(b_1,\ldots,b_r)\in (\Zm)^r$ is $\sum_{i=1}^r |a_i-b_i|$,
where $|x|=\min(x,m-x)$ for $x\in \Zm$, so it coincides with
$\|(a_1-b_1,\ldots,a_r-b_r)\|_2$, where $\|\cdot\|_2$ is as defined above.) The
{\em covering radius} is the smallest $\rho$ such that every vector of
$(\Zm)^r$ is $\rho$-covered. 

Now let ${\bf e}$ be the all one vector of $(\Zm)^r$. Obviously, for the
covering radius~$\rho$ of the code $C=(\Zm){\bf e}$ in the Lee metric we have
$$
  \rho=\begin{cases} g(m,r) & \text{ if } m=2, \\
                     h(m,r) & \text{ if } m>2.
       \end{cases}
$$
The Lee distance, and hence the covering radius based on it, is in general 
different from the Hamming distance; they coincide when $m=2$.

We can also interprete $g(m,r)$ and $h(m,r)$ as diameters of the graphs with
vertex set $V'$ where two vertices $\alpha$ and $\beta$ are connected if and
only if $\alpha-\beta\in S$ or $\in S\cup -S$, respectively (cf.\ \cite{gaso}
for prime $m$). Here $S$ is the set of generators $\{g_1,\ldots,g_r\}$ of
$V'$.

\section{Proof of Theorem \ref{ThmG}} \label{SecG}

We must solve the linear program that asks to maximise $\|{\bf v}\|_1$ under
$$
  \|{\bf v}\|_1 \le \|{\bf v}+x\cdot {\bf e}\|_1 \quad \text{ for all }
  x\in\Zm.
$$
Now since
$$
  \| {\bf v} + x{\bf e} \|_1 \equiv \| {\bf v} \|_1 + r\bar{x} \pmod m,
$$
the conditions of the linear program may be sharpened to
\begin{equation}
\label{tosum}
  \| {\bf v} \|_1 \le \| {\bf v} + x{\bf e} \|_1 - \overline{rx}
        \quad \text{ for all } x\in \Zm,
\end{equation}
where $\overline{rx}$, as above, denotes the remainder of $rx$ upon division by
$m$. Since each coordinate of ${\bf v}+x{\bf e}$ runs through all elements of
$\Zm$ as $x$ runs through $\Zm$, summing \eqref{tosum} over $x\in \Zm$ yields
\begin{align*}
  m\|{\bf v}\|_1 
           &\le r\sum_{x\in \Zm} \overline{x}-\sum_{x\in \Zm}\overline{rx}\\
           &= \binom{m}{2}r-\gcd(m,r)^2 \binom{m/\gcd(m,r)}{2}\\
           &= \frac{m((r-1)(m-1)+\gcd(m,r)-1)}{2}.
\end{align*}
Obviously, this upper bound is attained by a vector ${\bf v}$ with
$$
  t_k =_{\rm def} \frac{\overline{r(k-1)}+r-\overline{rk}}{m}$$
coordinates equal to $m-k$ for $k\in \Zm\setminus\{0\}$ and all other 
coordinates equal to zero. Namely, we have for $x\in \Zm\setminus\{0\}$, 
\begin{align*}
  \|{\bf v}+x{\bf e}\|_1 &= \|{\bf v}+(x-1){\bf e}\|_1 + r - mt_x \\
         &= \|{\bf v}\|_1 + \overline{r(x-1)} + r - mt_x\\
         &= \|{\bf v}\|_1 + \overline{rx} 
\end{align*}
by induction and thus equality in $(\ref{tosum})$.
\proofend

\section{Upper bounds} \label{SecUpper}

In this section and the next we prove Theorem \ref{ThmMR}.
Propositions \ref{PropRBig} and \ref{PropRSmall} will show that the values
taken by the function $h(m,r)$ indeed give an upper bound for the norm
$\|\cdot\|_2$ of an admissible vector in all cases. 
Throughout this section, we will write $\|\cdot\|$ for
$\|\cdot\|_2$.

We start with some preliminary results.

\begin{lemma} \label{LemSum}
We have
${\displaystyle
  \sum_{x\in\Z/m\Z} |x| = \begin{cases}
     \frac{m^2}4 & \text{ if $m$ is even} \\
     \frac{m^2-1}4 & \text{ if $m$ is odd}.
     \end{cases}
 }
$
\end{lemma}
The proof is left to the reader.

\begin{lemma} \label{LemCongruence}
Let $m$ be even. Then for any ${\bf v}\in V$, we have 
$\|\, {\bf v}+x\cdot {\bf e} \,\| \equiv \|{\bf v}\| + rx \pmod{2}$ 
for all $x\in \Z/m\Z$.
\end{lemma}

\paragraph*{Proof.} 
For even $m$, we have $|c+x| \equiv |c|+x\pmod{2}$ for all $c,x\in
\Z/m\Z$.  \proofend

~

The following Proposition gives upper bounds that are the right ones whenever
$r\ge m$, and also whenever $r$ is even. For the cases where $r$ is odd and
less than $m$, the bounds given in Proposition \ref{PropRSmall} are better
(see also Section \ref{SecProof}).

\begin{prop} \label{PropRBig}
Let ${\bf v}\in V$ be admissible. We have
$$
  \|{\bf v}\| \le \begin{cases}
    \frac{mr}4 - \frac{r}{4m} & \text{ if $m$ is odd} \\
    \frac{mr}4          & \text{ if $m$ and $r$ are both even} \\
    \frac{mr}4 - \frac12 & \text{ if $m$ is even and $r$ is odd}.
    \end{cases}
$$
\end{prop}

\paragraph*{Proof.}
The inequalities $\|{\bf v}+x{\bf e}\| \ge \|{\bf v}\|$ are summed over all
$x\in\Z/m\Z$.  By Lemma \ref{LemSum}, this yields
$$
  m\|{\bf v}\| \le \begin{cases} \tfrac{(m^2-1)r}4 & \text{ if $m$ is odd} \\
                                 \tfrac{m^2r}4     & \text{ if $m$ is even}.
             \end{cases}
$$
This can be sharpened if $m$ is even and $r$ is odd. Namely, by Lemma 
\ref{LemCongruence}, we find the sharper inequality
$$ 
  \|{\bf v}+x{\bf e}\| \ge \|{\bf v}\| + (x \!\!\!\mod 2);
$$
by summing over $x$, we get
$$
  m\|{\bf v}\| \le \tfrac{m^2r}4 - \tfrac{m}2.
$$
Now division by $m$ yields the result in all cases. \proofend

~

We now embark on the subcase where the dimension $r$ is odd and at most equal
to $m$, as we will need to strengthen the bounds in Proposition \ref{PropRBig}
for this case. Here, much more preparation is needed; the argument is concluded
in Proposition \ref{PropRSmall}.

\begin{defi} 
For a vector ${\bf v}\in V$, we define the \emph{norm sequence} of ${\bf v}$,
written $(N_x({\bf v}))$ or simply $(N_x)$ where $x$ runs over $\Z/m\Z$,
by setting $N_x = \|{\bf v}+x{\bf e}\|$. 
\end{defi}

\begin{lemma} \label{LemDistinct}
Let $r$ be odd, and let ${\bf v}\in V$. If $m$ is even, then $N_{x+1} \ne N_x$
for all $x\in\Z/m\Z$. If $m$ is odd and the number of distinct components of
${\bf v}$ is $s$, then there are at most $s$ values of $x$ in $\Z/m\Z$ for
which $N_{x+1}=N_x$. 
\end{lemma}

\paragraph*{Proof.}
For $m$ even, the result follows easily from Lemma \ref{LemCongruence}.

Suppose $m$ is odd. As $r$ is odd, we cannot have $N_{x+1}=N_x$ unless we have
$|v_i+x+1|=|v_i+x|$ for at least one $i$ with $1\le i \le r$. But this implies
$v_i+x=\frac{m-1}2$. Therefore, if ${\bf v}$ has $s$ distinct components, there
can exist at most $s$ distinct $x\in\Z/m\Z$ with $N_{x+1}=N_x$.
\proofend

~

The next two Lemmas deal with the horizontal symmetry or near-symmetry of the
norm sequence; they are applied in Lemma \ref{LemNormStrip}. The detailed first
assertions of both are again used in Section \ref{SecEvenDim}. For
$x\in\Z/m\Z$,  we will write $\bar{x}$ for the representative of $x$ in the set
$\{0,1,\ldots,m-1\}\subseteq \Z$, as before.

\begin{lemma} \label{LemMirrorEven}
Let $m$ be even. For all $x\in \Z/m\Z$, we have 
$$
  |x| + |x+\tfrac{m}{2}| = \tfrac{m}{2}.
$$ 
For all ${\bf v}\in V$, we have 
$$
  \|{\bf v}\| + \|{\bf v}+ \tfrac{m}2\cdot {\bf e} \| = \tfrac{mr}2.
$$
\end{lemma}

\paragraph*{Proof.}
If $0\le \bar{x} < \frac{m}{2}$, then $|x| + |x+\frac{m}{2}| = \bar{x} + m-
(\bar{x}+\frac{m}{2}) = \frac{m}{2}$.  If $\frac{m}{2} \le \bar{x} < m$, then
$|x| + |x+\frac{m}{2}| = m-\bar{x} + (\bar{x}+\frac{m}{2}-m) = \frac{m}{2}$.
The last assertion follows by the definition of the Lee norm. 
\proofend

\begin{lemma} \label{LemMirrorOdd}
Let $m$ be odd. For all $x\in \Z/m\Z$, we have 
\begin{romanlist}
  \item $|x| + |x+\frac{m+1}{2}| = 
        \begin{cases}
          \frac{m-1}2 & \text{if} \quad 0\le \bar{x} \le \frac{m-1}2 \\
	  \frac{m+1}2 & \text{if} \quad \frac{m+1}2 \le \bar{x} \le m-1;
        \end{cases}$
  \item $|x| + |x+\frac{m-1}{2}| = 
        \begin{cases}
          \frac{m-1}2 & \text{if} \quad x=0 \,\text{ or }\,
	                          \frac{m+1}2 \le \bar{x} \le m-1 \\
	  \frac{m+1}2 & \text{if} \quad 1 \le \bar{x} \le \frac{m-1}2.
        \end{cases}$
  \end{romanlist}
For all ${\bf v}\in V$, we have 
$$
  2\|{\bf v}\| + \left\|{\bf v}+\tfrac{m-1}2 \cdot {\bf e}\right\| + 
  \left\|{\bf v}+ \tfrac{m+1}2 \cdot {\bf e}\right\| 
  = mr - \#\{i \mid v_i=0 \}.
$$
\end{lemma}

\paragraph*{Proof.}
If $0\le \bar{x}\le \frac{m-1}2$, then $|x|=\bar{x}$ and $|\bar{x} +
\frac{m+1}2 | = m - (\bar{x}~+~\frac{m+1}2)$; if $\frac{m+1}2 \le \bar{x} \le
m-1$, then $|x|=m-\bar{x}$ and $|x + \frac{m+1}2 | = (\bar{x}~+~\frac{m+1}2
)-m$.

We have $|0| + |\frac{m-1}2 | = \frac{m-1}2 $. Also, if $1 \le \bar{x} \le
\frac{m-1}2 $, then $|x|=\bar{x}$ and $|x + \frac{m-1}2 | = m -
(\bar{x}~+~\frac{m-1}2 )$.  Finally, if $\frac{m+1}2 \le \bar{x} \le m-1$, then
$|x|=m-\bar{x}$ and $|x + \frac{m-1}2 | = (\bar{x}~+~\frac{m-1}2 )-m$.

As to the last assertion, let ${\bf v}=(v_1,\ldots,v_r)\in V$ and let $1\le i
\le r$. By the first part, we have
$$
  \left( |v_i| + |v_i + \tfrac{m-1}2 | \right) +
  \left( |v_i| + |v_i + \tfrac{m+1}2 | \right) = m,
$$
unless the two summands are equal. Now these two summands being both equal to
$\frac{m-1}2 $ implies $v_i=0$, and they cannot be both equal to $\frac{m+1}2
$. The claim follows by the definition of the Lee norm. 
\proofend

\begin{lemma} \label{LemNormStrip}
Let ${\bf v}\in V$ be admissible. Then for all $x\in \Z/m\Z$, we have
$$
  \| {\bf v} + x{\bf e} \| \le \tfrac{mr}2 - \|{\bf v}\|.
$$
\end{lemma}

\paragraph*{Proof.} 
First, suppose that $m$ is even, and apply Lemma \ref{LemMirrorEven} to
${\bf v}+x{\bf e}$. By admissibility, we have $\|{\bf v}+(x+\frac{m}2){\bf e}\|
\ge \|{\bf v}\|$, and the result follows.

If $m$ is odd, we apply Lemma \ref{LemMirrorOdd} to ${\bf v}+x{\bf e}$ and use
the admissibility inequality for both ${\bf v}+(x+\frac{m-1}2 ){\bf e}$ and
${\bf v}+(x+\frac{m+1}2){\bf e}$. After dividing by $2$, we obtain the result.
\proofend

\begin{defi} \label{DefExtr}
Let $(a_x)_{x\in\Z/m\Z}$ be a sequence of real numbers. We define the
\emph{slope} of $(a_x)$ at $x$ to be $a_{x+1}-a_x$. We say that the sequence
has a \emph{maximum} at $x$ if there exists $c\in \{1,2,\ldots,m-1\}$ such that
\begin{align*}
  a_{x-1} &< a_x; \\
  a_{x+i} &= a_x \text{ for } i=0,1,\ldots,c-1; \\
  a_{x+c} &< a_x.
\end{align*}
A \emph{minimum} is defined symmetrically; and we define an \emph{extremal
value} to be either a minimum or a maximum.
\end{defi}

\begin{lemma} \label{LemBend}
Let ${\bf v}\in V$, and let $(N_x)$ be the norm sequence of ${\bf v}$. If the
number of distinct components of ${\bf v}$ is $s$, then the number of extremal
values of the sequence $(N_x)$ is at most $2s$.
\end{lemma}

Note that this result is independent of the parities of $m$ and $r$. For the
multiplication by $2$ used in the proof of the second part, see also
Section \ref{SecOddOdd}.

\paragraph*{Proof.}
Recall that all sequences in this proof are periodic with period $m$.
The sequence $(N_x)$ is the sum of the sequences $(|v_i+x|)$, where $i$ runs
over $1,\ldots,r$. 

First, let us consider the case where $m$ is even. Here each period of the
composing sequences is made up of two segments; in the first, starting at
$x=\overline{-v_i}$, the sequence increases with slope $1$, while in the second
it decreases with slope $-1$. We see that the composing sequences only change
slope at the two extremal values they possess, which all have $c=1$ in the
notation of Definition \ref{DefExtr}. Now suppose $(N_x)$ has an extremal value
at $x$; then in particular its slope at $x-1$ and its slope at $x$ are
different, so one of the composing sequences must change its slope as well. It
follows that also one of the composing sequences has an extremal value at $x$,
and consequently $x$ must be equal to one of the at most $2s$ values where such
an extremal value occurs.

Second, assume $m$ is odd; we will reduce this case to the previous one, as
follows. Let $(S_x)$ be any sequence of real numbers indexed by the integers
modulo $m$, and suppose $(T_y)$ is any real sequence, indexed by the integers
modulo $2m$, such that $T_y=S_{y/2}$ whenever $y$ represents an even class
modulo $2m$. We claim that the sequence $(T_y)$ has no fewer extremal values
than the sequence $(S_x)$. Indeed, suppose $(S_x)$ has a maximum at $x$, and
consider the subsequence $T_{2x-2}=S_{x-1}, T_{2x-1}, T_{2x}=S_x, \ldots,
T_{2x+2c-1}, T_{2x+2c}$ of $(T_y)$. Let $y$ be the first index with $T_y$ as
large as possible in this subsequence. Then as $T_{y-1}<T_y$ and
$T_{2x+2c}<T_y$, the sequence $(T_y)$ has a maximum at $y$, possibly with
a smaller value of $c$. This proves the claim.

We apply the claim to the norm sequence $(N_x)$ of ${\bf v}$ and the sequence
$(M_y)_{y\in \Z/2m\Z}$ with $M_y=\frac12 \|2{\bf v}+y{\bf e}\|$ for
$y\in\Z/2m\Z$; here $2{\bf v}$ means the image of ${\bf v}$ under the
$\Z$-linear map $(\Z/m\Z)^r\rightarrow (\Z/2m\Z)^r$ that in every coordinate
maps $z$ to $2z$, for all $z\in\Z/m\Z$. Note that the norms $(M_y)$ are
evaluated modulo $2m$, whereas the $(N_x)$ are evaluated modulo $m$.  Clearly,
we have $N_x=M_{2x}$ for all $x\in \Z/m\Z$, so the claim applies. By the first
part, the sequence $(M_y)$ has at most $2s$ extremal values; consequently, the
same holds for the norm sequence $(N_x)$ of ${\bf v}$, and the Lemma is proved.
\proofend

~

We are now in a position to prove the upper bounds from Theorem \ref{ThmMR} in
the case where $r\le m$ and $r$ is odd.

\begin{prop}  \label{PropRSmall}
Let $r$ be odd, assume $r\le m$, and let ${\bf v}\in V$ be admissible.
Then we have
$$
  \|{\bf v}\| \le \frac{mr}4 - \frac{m}{4r}.
$$
\end{prop}

\paragraph*{Proof.}
Consider the norm sequence $(N_x)_{x\in \Z/m\Z}$ of ${\bf v}$. By leaving out
all members $x$ of the index set that have $N_{x-1}=N_x$, we arrive at a
subsequence $(N'_y)_{y\in\Z/m'\Z}$ of $(N_x)$, with period $m'\le m$. Note that
we no longer have $N'_y=\|{\bf v}+y{\bf e}\|$, because the $N'_y$ have been
renumbered. The subsequence has the following properties:
\begin{list}{\makebox[\labelwidth][r]{{\rm (\roman{romancnt})}}}
               {\usecounter{romancnt}%
                \topsep         5pt%
                \parskip        0pt%
                \partopsep      5pt%
                \itemsep        5pt%
                \parsep         0pt}
  \item $N'_y$ is a nonnegative integer for all $y$;
  \item we have $N'_{y+1} \ne N'_y$ for all $y$;
  \item the period $m'$ is equal to $m$ if $m$ is even, and is at least $m-r$
        otherwise;
  \item we have $\|{\bf v}\| \le N'_y \le \lfloor \frac{mr}2 \rfloor - \|{\bf
        v}\|$ for all $y$;
  \item the sequence $(N'_y)$ has at most $2r$ extremal values.
\end{list}
The last three of these follow by Lemmas \ref{LemDistinct}, \ref{LemNormStrip},
and \ref{LemBend}.

Now it is easy to see that if a sequence of integers is squeezed between
bounds $B$ from above and $A$ from below and cannot repeat itself, it must
have an extremal value at least every $B-A$ elements. Therefore, the number of
extremal values times the ``band width'' $B-A$ provides an upper bound on the
length of such a sequence. (With a finite sequence, there are some caveats at
the end points, but our sequences are periodic, and hence do not have end
points.)

We find therefore
\begin{align*}
  (2r)\left( \tfrac{mr}2 - 2\|{\bf v}\| \right) & \ge m'= m && 
    \text{ if $m$ is even, and} \\
  (2r)\left( \tfrac{mr}2 - \tfrac12 - 2\|{\bf v}\| \right) &\ge m'\ge m - r && 
    \text{ if $m$ is odd}.
\end{align*}
It turns out that the inequalities for the two cases are equivalent. The
result follows easily. \proofend

~

Note that the argument could be adapted to yield an upper bound also in the
cases where $r>m$. However, the resulting bound $\|{\bf v}\|\le \frac{mr}4 -
\frac14$ is larger than the ones given by Proposition \ref{PropRBig}. For
$m=r$, the two bounds coincide.

\section{Constructions} \label{SecCon}

After having shown that the values taken by the norm bound function $h(m,r)$
are upper bounds for the norms of admissible vectors, we will now proceed to
construct admissible vectors for all $m$ and $r$, the norm of which actually
attains these values. As in the last section, we write $\|\cdot\|$ for the
function $\|\cdot\|_2$, as defined in Section \ref{SecComb}.

\subsection{Even dimension} \label{SecEvenDim}

The case where the dimension $r$ is even, is relatively easy. In this case,
a useful building block for admissible vectors of high norm is the
\emph{optimal pair}. To achieve flexibility in constructions, we do not require
that an optimal pair be itself admissible.

\begin{defi}
An \emph{optimal pair} is a vector ${\bf v}$ of length $2$ such that for some
$x\in\Z/m\Z$, the vector ${\bf v}+x{\bf e}$ is admissible of maximal norm.
\end{defi}

\begin{lemma} \label{LemEven2}
If $m$ is even, then for all $y\in\Z/m\Z$, the vector $(y,y+\frac{m}{2})$ is an
optimal pair, and is admissible of norm $\frac{m}{2}$.
\end{lemma}

\paragraph*{Proof.}
For all $x\in \Z/m\Z$, we have $\| (y,y+\frac{m}{2}) + (x,x) \| = |y+x| +
|y+x+\frac{m}{2}| = \frac{m}{2}$, by Lemma \ref{LemMirrorEven}. This norm
is maximal by Proposition \ref{PropRBig}.
\proofend

\begin{lemma} \label{LemOdd2}
If $m$ is odd, then for all $y\in \Z/m\Z$ the vector $(y,y+\frac{m-1}2)$ is
an optimal pair. When $y=0$ or $\frac{m+1}2\le \bar{y} \le m-1$, such a vector
is admissible of norm $\frac{m-1}2$.
\end{lemma}

\paragraph*{Proof.}
The assertions follow directly from Lemma \ref{LemMirrorOdd}, with Proposition
\ref{PropRBig} showing that the attained norm is maximal.
\proofend

\paragraph*{}
The next result shows that the bounds of Proposition \ref{PropRBig} are
sharp in the case that the dimension $r$ is even. 

\begin{prop} \label{PropEven}
Let $r$ be even.
\begin{romanlist}
  \item If $m$ is even, then there exists an admissible vector ${\bf v}$ of
        length $r$ and norm $\frac{mr}{4}$.
  \item If $m$ is odd, then there exists an admissible vector ${\bf v}$ of
	length $r$ and norm $\left\lfloor \frac{mr}4 - \frac{r}{4m}
	\right\rfloor$.
\end{romanlist}
\end{prop}

\paragraph*{Proof.} 
For even $m$, the vector 
${\bf v}=\left(0,\frac{m}{2}\right)^{r\!/2}=\left(0,\frac{m}{2},0,\frac{m}{2},\ldots,0,\frac{m}{2}\right)$ 
is clearly
admissible of the given norm, by Lemma \ref{LemEven2} and the fact that the
concatenation of admissible vectors yields again an admissible vector.

For the case of odd $m$, we use Lemma \ref{LemOdd2} and the same fact, with
some subtility. Let ${\bf v}=(y,y+ \frac{m-1}2)$ be an optimal pair for $m$,
and let $(N_x(v))_{x\in\Z/m\Z}$ be its norm sequence. From Lemma
\ref{LemMirrorOdd}, it is easy to see that we have
$$
  N_x({\bf v}) = \begin{cases}
             \frac{m+1}2 & \text{ if } 
	             x\in\{ -y+1,-y+2,\ldots,-y+\frac{m-1}2\}; \\
             \frac{m-1}2 & \text{ if } 
	             x\in\{ -y + \frac{m+1}2,\ldots,-y + m\}.
          \end{cases}
$$
We will call these two subsets of $\Z/m\Z$ the \emph{high} and \emph{low
regions} of $N_x({\bf v})$, respectively.

We will determine $r/2$ optimal pairs such that their concatenation is
admissible of maximal norm. For this, it is necessary to select the pairs in
such a way that the high regions of their norm sequences are spread as evenly
as possible over the total range $x=0,\ldots,m-1$.

Writing ${\bf v}_i = (y_i, y_i + \frac{m-1}2)$, we take $y_i=-(i-1)\frac{m-1}2 $ for $i\ge 1$. The high region of $(N_x({\bf v}_i))$ 
starts at
$x=(i-1)\frac{m-1}{2}+1$
and ends at 
$x=i\frac{m-1}{2}$.
We see that the high regions of $r/2$ pairs, put in a row, cover a contiguous
region from $x=1$ to $x=\frac{r}2 \frac{m-1}2$; reducing the indices modulo
$m$, we find that every element in the range $x=0,\ldots,m-1$ is covered at
least
$$
  \left\lfloor \frac{\frac{r}2 \frac{m-1}2 }{m} \right\rfloor
$$
times. Moreover, at $x=0$, and possibly some elements to the left of $x=0$,
this inequality is an equality, because covering ``started'' at $x=1$, strictly
to the right of $x=0$. This means that the concatenation ${\bf v}$ of the pairs ${\bf v}_i$
thus selected is admissible, and that its norm satisfies
$$
  \|{\bf v}\| \ge \tfrac{r}2 \tfrac{m-1}2 + 
    \left\lfloor \frac{\frac{r}2 \tfrac{m-1}2 }{m} \right\rfloor
        = \left\lfloor \tfrac{r(m-1)}4 + \tfrac{r(m-1)}{4m} \right\rfloor
        = \left\lfloor \tfrac{mr}4 - \tfrac{r}{4m} \right\rfloor.
$$
By Proposition \ref{PropRBig}, we must have equality here, and the construction
is finished.~\proofend

\subsection{Odd dimension, even modulus}

We now proceed to the case of odd dimension, which is more complicated. We
first assume that $m$ is even, and that $r<2m$. The construction of an
admissible vector for such parameters is derived from the proof of Proposition
\ref{PropRSmall}; we try to choose the components of a vector ${\bf
v}=(v_1,\ldots, v_r)$ such that its norm sequence $(N_x({\bf v}))$ has always
slope $\pm 1$ and has its extremal values spread as evenly as possible over the
range $x=0,\ldots, m-1$. As earlier, we write $V=(\Z/m\Z)^r$, and for $x\in
\Z/m\Z$, we write $\bar{x}$ for the representative of $x$ in the set
$\{0,\ldots,m-1\}\subseteq \Z$.

\begin{defi}
  Assume $m$ even and $r$ odd. A vector ${\bf v}\in V$ satisfying
  \begin{equation} \label{EqDistr} 
    0=\bar{v}_1\le \bar{v}_2-\tfrac{m}2 \le \bar{v}_3 \le \bar{v}_4-\tfrac{m}2
    \le \ldots \le \bar{v}_{r-1}-\tfrac{m}2 \le \bar{v}_r < \tfrac{m}2
  \end{equation}
  will be called \emph{balanced}.
\end{defi}

\begin{lemma} \label{LemDistr}
  Let ${\bf v}\in V$ be balanced and let $(N_x)$ be its norm sequence. Then
  we have $N_{x+1}-N_x=\pm 1$ for all $x\in\Z/m\Z$.
\end{lemma}

\paragraph*{Proof.}
  As $m$ is even, each individual component $v_i$ has $|v_i+x+1|-|v_i+x|=
  \pm 1$ for all $x$, the sign being positive when $v_i+x=0,1,\ldots,\tfrac{m}2
  -1$ and negative otherwise. At $x=0$, we have exactly $\tfrac{r+1}2$
  ``increasing'' and $\tfrac{r-1}2$ ``decreasing'' components, so that
  $N_1-N_0=1$.

  But by the alternating arrangement of the $v_i$ around $\tfrac{m}2$, it is
  clear that after a component changes from increasing to decreasing at a
  certain $x$, we cannot have another component doing the same; we must first
  see a component changing from decreasing to increasing, possibly at the same
  $x$ if the corresponding inequality in \eqref{EqDistr} is an equality. Thus,
  the balance between increasing and decreasing components is always either $1$
  or $-1$, and the assertion is clear.~\proofend

~

We have shown earlier (Lemma \ref{LemBend}) that the norm sequence of any
vector ${\bf v}$ in $V$ has at most $2s$ extremal values, where $s$ is the
number of distinct components of ${\bf v}$. Now assume ${\bf v}$ is balanced.
Then in fact, an extremal value will occur whenever the balance between the
numbers of increasing and decreasing components of ${\bf v}$ changes. For this,
we look at the extremal values of the composing sequences. If $i$ is odd, then
$0\le \bar{v}_i < \frac{m}2$, so the sequence $|v_i+x|$ has a maximum at
$x=\frac{m}2-v_i$. If $i$ is even, then $\frac{m}2 \le \bar{v}_i < m$, so a
minimum occurs at $x=m-v_i$.  All these values for $x$ are possible locations
of extremal values in the norm sequence of ${\bf v}$. Counting from $x=1$
onwards, the first location is $\frac{m}2-v_r$, the second is $m-v_{r-1}$, and
so on.  Finally, we start by having a minimum at $x=0$.

Thus, let us define
\begin{equation} \label{EqM}
  \begin{aligned}
    m_0 &= \|{\bf v}\|; \\
    m_i &= \begin{cases} \left\|{\bf v} + 
                           \left(\tfrac{m}2-v_{r-i+1}\right){\bf e}\right\| 
                           & \text{ if $i$ is odd}; \\
			 \left\|{\bf v} + 
			   \left(m-v_{r-i+1}\right){\bf e}\right\|
                           & \text{ if $i$ is even}.
  	   \end{cases}
  \end{aligned}
\end{equation}
Then the $m_i$, for $i=0,\ldots,r$, include all extremal values of the norm
sequence $(N_x)$ of ${\bf v}$ in the range $x=0,\ldots,\frac{m}2$.

\begin{lemma} \label{LemMDiff}
  Let ${\bf v}$ be balanced. Then we have $m_1-m_0=\tfrac{m}2 - \bar{v}_r$,
  while for $i=1,\ldots,r-1$, 
  $$
    m_{i+1} - m_i = \begin{cases}
		       (\bar{v}_{r-i} - \tfrac{m}2) - \bar{v}_{r-i+1} 
		           & \text{ if $i$ is odd}; \\
		       (\bar{v}_{r-i+1} - \tfrac{m}2) - \bar{v}_{r-i} 
		           & \text{ if $i$ is even}.
		    \end{cases}
  $$
\end{lemma}

\paragraph*{Proof.}
  First assume $i$ is odd; then $m_i$ is a possible maximum of the norm
  sequence, occurring at $x=\frac{m}2-v_{r-i+1}$. The subsequent possible
  minimum $m_{i+1}$ occurs at $x=m-v_{r-i}$. If these values for $x$ are
  equal, then we also have $m_{i+1}=m_i$ and the claim is proved. If not,
  then between these values of $x$ the norm sequence has a constant slope
  of $-1$ (cf.\ Lemma \ref{LemDistr}). Therefore, the difference 
  $m_{i+1}-m_i$, as claimed, is equal to
  $$
    (-1) \cdot \left( ( m-\bar{v}_{r-i} ) - \left( \tfrac{m}2 - \bar{v}_{r-i+1}
    \right) \right).
  $$
  The case where $i>0$ is even and the case $i=0$ are analogous.
\proofend

\begin{lemma} \label{LemBand}
  Let ${\bf v}$ be balanced. Then ${\bf v}$ is admissible if and only if
  $N_0\le m_i\le m_r=\frac{mr}2-N_0$ for all $i$.
\end{lemma}

\paragraph*{Proof.}
  We continue to assume $m$ even and $r$ odd; by definition, we have $m_r=
  N_{m/2}$. Now we use the symmetry in the norm sequence given by Lemma
  \ref{LemMirrorEven}, which says that, for all $x$,
  $$
    N_{x+\tfrac{m}2} = \tfrac{mr}2 - N_x.
  $$
  First assume ${\bf v}$ is admissible; then from $N_0\le N_{x+\tfrac{m}2}$, we
  find $N_x\le N_{m/2}$ by using the formula twice. Thus in particular all
  $m_i$ are between $N_0$ and $m_r=N_{m/2}=\frac{mr}2-N_0$, as claimed.

  For the other direction, from $N_0\le m_i\le m_r$ for all $i$, we find
  $N_0\le N_x\le m_r$ for $x\le m/2$, because the $m_i$ contain among
  themselves all extreme values of the first half of the sequence $(N_x)$. But
  then by symmetry $N_{x+\tfrac{m}2} = \tfrac{mr}2 - N_x \ge \tfrac{mr}2 - m_r
  =N_0$, so we have $N_0\le N_x$ for all $x$, as desired.
\proofend

~

The next Lemma shows that there are several equivalent options for the
formulation of the norm bound function in \eqref{EqDefH}, when $m$ is even and
$r$ odd, and $r$ is not too far away from $m$. In fact, comparable formulae can
be given in case $m$ is odd also, but we omit these as they are not needed in
the sequel. The proof is left to the reader.

\begin{lemma} \label{LemProperties}
  Let $m$ be even and $r$ odd, and assume $\frac{m}2 \le r \le 2m$. Then
  $$
    \lfloor \tfrac{mr}4 - \tfrac{m}{4r} \rfloor = \lfloor \tfrac{mr}4
    - \tfrac12 \rfloor = \lfloor \tfrac{mr}4 - \tfrac{r}{4m} \rfloor =
    \begin{cases}
      \tfrac{mr}4 - 1 & \text{ if } m \equiv 0 \pmod{4}; \\
      \tfrac{mr}4 - \tfrac12 & \text{ if } m \equiv 2 \pmod{4}.
    \end{cases}
  $$
\end{lemma}

\begin{lemma} \label{LemC}
  Let $m$ be even and $r$ odd, with $r \le 2m$, and let $Q$ and $R$ be integers
  such that 
  $$
    \tfrac{m}2 = Qr + R, \text{ with } 0\le R<r.
  $$
  Then the quantity $C=\tfrac{mr}2 - 2h(m,r)$ satisfies
  $$
    C=\begin{cases}
        Q    & \text{ if } R=0; \\
	Q+1  & \text{ if $R$ is odd}; \\
	Q+2  & \text{ if $R$ is positive and even}.
      \end{cases}
  $$
  Furthermore, we have $C\equiv \frac{m}2 \pmod{2}$.
\end{lemma}

\paragraph*{Proof.} 
  Recall that $h(m,r) = \lfloor \frac{mr}4 - \frac{m}{4r} \rfloor$ with our
  assumptions, by \eqref{EqDefH} and Lemma \ref{LemProperties}. The proof is
  tedious but easy, and is left to the reader.
\proofend

\begin{prop} \label{PropMEvenROdd}
  Assume $m$ is even and $r$ is odd, with $r\le 2m$. Then there exists an
  admissible vector ${\bf v}\in V$ of norm $h(m,r)=\lfloor \frac{mr}4 -
  \frac{m}{4r} \rfloor$.
\end{prop}

\paragraph*{Proof.}
  We want to construct a balanced vector $v$ satisfying the requirements. Then
  by Lemma \ref{LemBand}, we must choose the components $v_i$ of ${\bf v}$ such
  that the associated quantities $m_i$ satisfy
  \begin{equation} \label{EqBand}
    m_0 = h(m,r) \le m_i \le m_r=\tfrac{mr}2 -h(m,r).
  \end{equation}
  Together with the constraints \eqref{EqDistr}, this is an integer programming
  problem in the variables $\bar{v}_1,\ldots,\bar{v}_r$. By Lemma
  \ref{LemMDiff}, the differences $m_{i+1}-m_i$ are, up to sign and in reverse
  order, the same as the differences $(\bar{v}_{i+1}-\frac{m}2)-\bar{v}_i$ and
  $\bar{v}_{i+1}-(\bar{v}_i-\tfrac{m}2)$ of the quantities figuring in
  \eqref{EqDistr}.  Thus it is enough to specify the values of the $m_i$, as
  both $m_0$ and $v_1=0$ are fixed.

  An easy but useful corollary of Lemma \ref{LemMDiff}, proved using 
  telescoping sums, is that
  \begin{equation} \label{EqM2}
    \sum_{i=0}^{r-1} |m_{i+1}-m_i|= \tfrac{m}2.
  \end{equation}
  Let us write $C$ for the difference $\frac{mr}2-2h(m,r)$ of the largest and
  the smallest $m_i$. By Lemma \ref{LemC}, $C$ is equal to or slightly larger
  than $\frac{m}{2r}$. This observation, together with \eqref{EqBand}, suggests
  that we take the $|m_{i+1}-m_i|$ all approximately equal to $\frac{m}{2r}$.
  The rest of the proof will give \emph{exact integer values} for the $m_i$
  so as to solve the integer programming problem for the $\bar{v}_i$. We note
  that, as $m_1-m_0=m/2-\bar{v}_r>0$ by \eqref{EqDistr}, we cannot put
  $m_1-m_0=0$.

  Let $Q$ and $R$ be integers satisfying
  $$
    \tfrac{m}2 = Qr + R, \text{ with } 0\le R < r.
  $$

  If $R=0$, the solution is easy, as we simply put 
  $$
    m_{i+1}-m_i=(-1)^i Q \quad\text{ for } i=0,\ldots, r-1.
  $$
  By Lemma \ref{LemC}, we have $C=Q$ in this case, so that
  \eqref{EqBand} is satisfied.

  If $R\ne 0$, we put
  $$
    m_{i+1}-m_i= \begin{cases}
                   (-1)^i (Q+1) & \text{ for } i=0 \text{ and }
		              i=r-R+1,\ldots,r-1; \\
                   (-1)^i Q & \text{ for } i=1,\ldots,r-R.
		 \end{cases}
  $$
  If $R$ is then odd, this implies that $m_i=m_0+1$ for all even $i$ with $2\le
  i\le r-R$, and $m_i=m_0$ for the other even $i$; furthermore, by Lemma
  \ref{LemC} we have $C=Q+1$, and in fact $m_r=m_0+C$, as the number of $i$
  with $|m_{i+1}-m_i|=Q$, which is $r-R$, is even. If $R$ is even and positive,
  we have $m_i=m_0+1$ for all even $i\ge 2$. In this case, by Lemma \ref{LemC}
  we have $C=Q+2$, and in fact we get $m_r=m_0+Q+2$, as the number of steps of
  size $Q$ is then odd.

  It follows that the integer programming problem defining the $\overline{v_i}$
  always has a solution, so that the existence of the required vector is
  proved.
\proofend

\subsection{Odd dimension, odd modulus} \label{SecOddOdd}

We continue to assume that $r$ is odd. We will now reduce the case of odd
modulus $m$ to the even case, using \emph{division by $2$}; this seems to be
the easiest way of extending the argument used in the proof of Proposition
\ref{PropMEvenROdd}. For $r\le m$, we achieve this reduction in Corollary
\ref{CorMROdd} below. The case $r>m$ will be dealt with in Section
\ref{SecLargeOdd}.

The group homomorphism $\Z/m\Z \rightarrow \Z/2m\Z$ sending $1$ to $2$ induces
a linear map $\mu_2:(\Z/m\Z)^r \rightarrow (\Z/2m\Z)^r$ that multiplies all
components by $2$. The image of $\mu_2$ consists of those vectors in
$(\Z/2m\Z)^r$ that have all their components even; we will call these
\emph{even vectors}. The map $\mu_2$ has an inverse on the set of even vectors
that we shall call \emph{division by $2$} and denote by ${\bf v}\mapsto {\bf
v}/2$.

Note that $\|\mu_2(v)\|$, as evaluated in $(\Z/2m\Z)^r$, is equal to $2\|v\|$,
when evaluated in $(\Z/m\Z)^r$, so that the Lee norm is multiplied by $2$ under
the map $\mu_2$; likewise, division by $2$ halves the norm.

\begin{lemma}
  If ${\bf v}\in (\Z/2m\Z)^r$ is even and admissible, then ${\bf v}/2\in
  (\Z/m\Z)^r$ is also admissible.
\end{lemma}

\paragraph*{Proof.}
  We have $\|{\bf v}\|\le \|{\bf v}+x\cdot {\bf e}\|$ for all $x\in\Z/2m\Z$; in
  particular, this holds for all \emph{even} $x\in \Z/2m\Z$, and so $\|{\bf
  v}/2\| \le \|{\bf v}/2 + x\cdot {\bf e}\|$ for all $x\in\Z/m\Z$.
\proofend

~

Recall that $h(m,r)$, as defined in \eqref{EqDefH}, gives the maximal norm of
an admissible vector of length $r$ and modulus $m$.

\begin{lemma} \label{LemR4}
  Let $m\equiv 2$ modulo $4$, and assume $r<2m$ and $r$ odd. Then
  $$
    h(m/2,r)=\left\lfloor \frac{h(m,r)}2 \right\rfloor.
  $$
  Furthermore, if $Q$ and $R$ are integers such that $\frac{m}2 = Qr+R$ with
  $0\le R<r$, then $h(m,r)$ is even if $R=0$ or $R\equiv 2\pmod{4}$ or $R\equiv
  r\pmod{4}$, and odd otherwise.
\end{lemma}

\paragraph*{Proof.}
  We use Lemma \ref{LemProperties} to have the formula $h(m,r)=\lfloor
  \tfrac{mr}4 - \tfrac{m}{4r} \rfloor$ from \eqref{EqDefH}, which holds
  for $r\le m$, also for $m < r < 2m$.

  Now we have $h(m/2,r)=\lfloor \frac{mr}8~-~\frac{m} {8r} \rfloor$ and
  $h(m,r)=\lfloor \frac{mr}4~-~\frac{m}{4r} \rfloor$. Because $\lfloor 
  \tfrac{x}2\rfloor = \lfloor \tfrac12 \lfloor x\rfloor\rfloor$ for any
  real $x\ge 0$, the first assertion easily follows.

  We now prove the second assertion. By substituting $2(Qr+R)$ for $m$ in the 
  formula for $h(m,r)$, we find
  $$
    h(m,r)= Q \cdot \frac{r^2-1}2 + 
            \begin{cases}
              0             & \text{ if } R=0; \\
              \tfrac{Rr-2}2 & \text{ if $R$ is nonzero and even}; \\
	      \tfrac{Rr-1}2 & \text{ if $R$ is odd}.
            \end{cases}
  $$
  The first term is even, so the parity of $h(m,r)$ equals the parity of the
  second term.
\proofend

\begin{prop}
  Let $m$ be congruent to $2$ modulo $4$, and assume $r \le m/2$ and $r$ odd.
  Then there exists in $V$ an \emph{even} admissible vector of norm
  $2h(m/2,r)$. 
\end{prop}

\paragraph*{Proof.}
  We will use the method developed in the proof of Proposition
  \ref{PropMEvenROdd} to construct a balanced admissible \emph{even} vector
  ${\bf v}$ satisfying the requirements.

  As above, we consider the components $v_i$ of ${\bf v}$ as the variables of
  an integer programming problem, which is here given by the constraints
  \eqref{EqDistr}, together with the following adaption of \eqref{EqBand}:
  \begin{equation} \label{EqBand2}
    m_0 = 2h(m/2,r) \le m_i \le m_r=\tfrac{mr}2 - 2h(m/2,r)
  \end{equation}
  for $i=0,\ldots,r$, and the additional constraint that all the $v_i$ must
  be even. Of course, as we fix $v_1=0$ and as $m/2$ is odd, this is
  equivalent to all the differences $(\bar{v}_{i+1}-\frac{m}2)-\bar{v}_i$ or
  $\bar{v}_{i+1}-(\bar{v}_i-\frac{m}2)$ being \emph{odd}, and this again to the
  differences $m_{i+1}-m_i$ being odd for all $i$ --- cf.\ Lemma
  \ref{LemMDiff}.

  Write $C'$ for the difference $m_r-m_0=\frac{mr}2 - 4h(m/2,r)$, and let $C$
  be as in Lemma \ref{LemC}. By Lemma \ref{LemR4}, we see that \eqref{EqBand}
  is equivalent to \eqref{EqBand2}, and we have $C'=C$, whenever $h(m,r)$ is
  even; if $h(m,r)$ is odd, this means that an even vector of norm $h(m,r)$
  does not exist, and we have to weaken \eqref{EqBand}, taking $C'=C+2$.
  
  As before, let $Q$ and $R$ be integers satisfying
  $$
    \tfrac{m}2 = Qr + R, \text{ with } 0\le R < r.
  $$
  We now have the same three cases, depending on whether $R$ is zero, odd, or
  nonzero and even. Again, we recall that we may not put $m_1-m_0=0$.
 
  First, suppose $R=0$. As $C=C'$ in this case, we have the same constraints as
  in the proof of Proposition \ref{PropMEvenROdd}. There, we gave
  $|m_{i+1}-m_i|$ the value $Q$ for all $i$. But $Q$ is odd, which means that
  we automatically obtain an even vector, and we are done.
 
  Now suppose $R$ is odd. We must distinguish two subcases. Thus, first suppose
  that $R$ and $r$ are congruent modulo $4$. It then follows by Lemma
  \ref{LemR4} that $C'=C=Q+1$. We cannot give $|m_{i+1}-m_i|$ the value $Q$
  now, as we did previously, since $Q$ is even. Instead, we take
  $$
    m_{i+1}-m_i = \begin{cases}
                     (-1)^i(Q+1) & \text{ for } i=0,\ldots,\tfrac{r+R}2-1; \\
                     (-1)^i(Q-1) & \text{ for } i=\tfrac{r+R}2,\ldots,r-1.
		  \end{cases}
  $$
  Note that by the assumption $r\le m/2$, we have $Q\ge 1$. Here we have
  $m_i=m_0$ for even $i\le \frac{r+R}2$ and $m_i=m_0+2$ for larger even $i$.

  If $R$ is odd, but not congruent to $r$ modulo $4$, we find by Lemma
  \ref{LemR4} that $h(m,r)$ is odd, and we have to take $C'=C+2=Q+3$. The
  assignment of values will be
  $$
    m_{i+1}-m_i = \begin{cases}
                     (-1)^i(Q+1) & \text{ for } i=0,\ldots,\tfrac{r+R}2-2; \\
		     (-1)^i(Q-1) & \text{ for } i=\tfrac{r+R}2-1,\ldots,r-2; \\
		     Q+3 & \text{ for } i=r-1.
		  \end{cases}
  $$
  Finally, suppose $R$ is nonzero and even. Again we find two subcases. Assume
  $R\equiv 2\pmod{4}$; then by Lemma \ref{LemR4} we find $C'=C=Q+2$. As we
  cannot assign the even value of $Q+1$, we take the assignment of values to be
  $$
    m_{i+1}-m_i = \begin{cases}
                    (-1)^i Q     & \text{ for } i=0,\ldots,r-\tfrac{R}2-1; \\
		    (-1)^i (Q+2) & \text{ for } i=r-\tfrac{R}2,\ldots,r-1.
		  \end{cases}
  $$
  The last case is where $R$ is nonzero and $R\equiv 0\pmod{4}$. By Lemma
  \ref{LemR4}, we see that $h(m,r)$ is odd and we must allow $C'=C+2=Q+4$ in
  \eqref{EqBand2} in order for an even vector to exist. Here, one can assign
  values of
  $$
    m_{i+1}-m_i = \begin{cases}
                    (-1)^i Q    & \text{ for } i=0,\ldots,r-2-\tfrac{R-4}2; \\
		    (-1)^i(Q+2) & \text{ for } i=r-1-\tfrac{R-4}2,\ldots,r-2;
		      \\
		    Q+4          & \text{ for } i=r-1.
		  \end{cases}
  $$
  In all the preceding cases, one checks easily that \eqref{EqBand2} is
  satisfied; the checks are the easier as we have chosen values for the $m_i$
  such that $m_i=m_0$ for all even $i$, except when $R\equiv r\pmod{4}$.
\proofend

\begin{cor} \label{CorMROdd}
  Let $m$ be odd, and assume $r \le m$ and $r$ odd. Then there exists in $V$
  an admissible vector of norm $h(m,r)$.
\end{cor}

\paragraph*{Proof.}
  Let ${\bf v}$ be an admissible even vector in $(\Z/2m\Z)^r$ of norm
  $2h(m,r)$, as provided by the Proposition; then ${\bf v}/2$ is the desired
  vector in $V$.
\proofend

\subsection{Large, odd dimension} \label{SecLargeOdd}

We just proved the norm bounds of Theorem \ref{ThmMR} sharp for $r$ odd and
at most equal to $2m$ (for $m$ even) or at most equal to $m$ (for $m$ odd). The
last step of the proof of the Theorem is to reduce the case of arbitrarily
large odd dimension to one of these cases, or to a case of even $r$. For this,
we use the fact that admissible vectors of maximal norm are particularly easy
to construct when the dimension $r$ is \emph{divisible} by the modulus $m$.

\begin{lemma} \label{LemFull}
Suppose $m$ divides $r$. Then the vector 
$$
  (0,1,\ldots,m-1)
$$
is admissible of maximal norm $\frac{m^2}{4}$ (if $m$ is even), resp.\
$\frac{m^2}4 - \frac14$ (if $m$ is odd).
\end{lemma}

\paragraph*{Proof.} 
Let ${\bf v}=(0,1,\ldots,m-1)$; adding ${\bf e}=(1,1,\ldots,1)$ to the vector
only permutes the coordinates, so it is clearly admissible. Its norm is given
by Lemma \ref{LemSum}.
\proofend

\begin{lemma} \label{LemReduce}
Suppose ${\bf v}$ is an admissible vector of length $r$ and maximal norm. If
$r\ge m$, then the concatenation of ${\bf v}$ with $(0,1,\ldots,m-1)$, of
length $r+m$, is also admissible of maximal norm. If $m$ is odd and $r$ is
even, this even holds for all $r\ge 1$.
\end{lemma}

\paragraph*{Proof.}
Write ${\bf w}$ for the concatenation of ${\bf v}$ with $(0,1,\ldots,m-1)$. We
use the fact that the concatenation of two admissible vectors is admissible,
with the norm of the concatenated vector being the sum of the norms of the two
summands. Therefore, it remains to prove that the concatenation again has
maximal norm.

According to Proposition \ref{PropRBig}, there are three cases. Now the
equalities
\begin{align*}
  \tfrac{mr}4 + \tfrac{m^2}4 &= \tfrac{m(r+m)}4, \\
  \left\lfloor \tfrac{mr}4 - \tfrac12 \right\rfloor + \tfrac{m^2}4 &=
  \left\lfloor \tfrac{m(r+m)}4 - \tfrac12 \right\rfloor, \text{ and} \\
  \left\lfloor \tfrac{mr}4 - \tfrac{r}{4m} \right\rfloor + \tfrac{m^2-1}4 &=
  \left\lfloor \tfrac{m(r+m)}4 - \tfrac{r+m}{4m} \right\rfloor
\end{align*}
settle the cases $m$ and $r$ both even, $m$ even and $r$ odd, and
$m$ odd, respectively.~\proofend

\begin{prop} \label{PropReduce}
  Let $m$ be given. If the norm bounds given in Theorem \ref{ThmMR} are sharp
  for $r$ with $1\le r \le 2m-1$, then they are sharp for all $r$.

  If the norm bounds are sharp for $m$ odd and $r$ even with $r\le m$, then
  they are also sharp for $r$ odd with $m<r\le 2m-1$.
\end{prop}

\paragraph*{Proof.}
Suppose we have $m$ and $r$ with $r\ge 2m$; write $r=Qm+R$ with integers $Q,R$
satisfying $m\le R < 2m$. An admissible vector of maximal norm of length $r$ is
constructed by concatenating such a vector of length $R$ with $Q$ copies
of $(0,1,\ldots,m-1)$, by Lemma \ref{LemReduce}.

As to the second statement, let $m$ and $r$ be odd with $m<r\le 2m-1$, and let
${\bf v}$ be an admissible vector of length $r-m$ and maximal norm. Then by the
last statement of Lemma \ref{LemReduce}, the concatenation of ${\bf v}$ with
$(0,1,\ldots,m-1)$ is admissible of length $r$ and maximal norm.
\proofend

\section{Proof of Theorems \ref{ThmMR}, \ref{Thm1}, and \ref{Thm2}} 
\label{SecProof}

\paragraph*{Proof of Theorem \ref{ThmMR}.}
Write $V=(\Z/m\Z)^r$, as before, and let $\|\cdot\|$ denote the \emph{norm} 
$\|\cdot\|_2$, as defined in Section \ref{SecComb}. We must prove that for all
$m$ and $r$, admissible vectors of norm $h(m,r)$ exist in $V$, and that
admissible vectors cannot have higher norms.

The fact that $h(m,r)$ forms an upper bound for the norm of an admissible
vector is proved in Propositions \ref{PropRBig}, for the cases where $r\ge m$
or $r$ is even, and \ref{PropRSmall} for the cases where $r$ is odd and $r\le
m$. In fact, if $r\le m$ and $m$ and $r$ both odd, it is clear that
$$
  \frac{mr}4 - \frac{r}{4m} \ge \frac{mr}4 - \frac{m}{4r};
$$
here the left hand side is the bound given by Proposition \ref{PropRBig}, and
the right hand is given by Proposition \ref{PropRSmall}. Also, if $m$ is even 
and $r$ odd, then for $r\le m/2$ the inequality
$$
  \frac{mr}4 - \frac12 \ge \frac{mr}4 - \frac{m}{4r}
$$
shows that the left bound, given by Proposition \ref{PropRBig}, is larger than
the right one from Proposition \ref{PropRSmall}, while for $m/2< r \le m$ the
floors of the two bounds are shown to be equal by Lemma \ref{LemProperties}.

The question whether the norm bound $h(m,r)$ is sharp was settled in Section 
\ref{SecCon}, in several cases, as follows.

For $r$ even, concrete vectors attaining the norm bound are given by
Proposition \ref{PropEven}. 

Assume $r$ is odd. By Proposition \ref{PropReduce}, we may reduce to a case
with $r<2m$, where the new $r$ can have either parity. Now if $r$ is even, we
use Proposition \ref{PropEven} to conclude the argument. If $r$ is odd and $m$
is even, we use Proposition \ref{PropMEvenROdd}. If both $m$ and $r$ are odd
and $m< r< 2m$, we use the second statement of Proposition \ref{PropReduce}
to conclude: the norm bound is sharp for modulus $m$ and even dimension $r-m$
by Proposition \ref{PropEven}, and hence it is sharp for modulus $m$ and odd
dimension $r$. If, finally, both $m$ and $r$ are odd and $r\le m$, we conclude
using Corollary \ref{CorMROdd}.
\proofend

~

We can now prove Theorems \ref{Thm1} and \ref{Thm2}.

\paragraph*{Proof of Theorem \ref{Thm1}.}
Note that the nonzero $(p^{r-1}-1)/r$th powers in $\F_{p^{r-1}}$ are exactly 
the $r$th roots of unity.

Now let $\xi$ be a primitive $r$th root of unity in $\F_{p^{r-1}}$. 
Since $p$ is a primitive root modulo~$r$, the field $\F_{p^{r-1}}$ is generated 
by $\xi$, i.\,e.\ $\{1,\xi,\ldots,\xi^{r-2}\}$ is a basis of $\F_{p^{r-1}}$ 
over $\F_p$. Since 
$$
  \sum_{i=0}^{r-1} \xi^i=0
$$ 
is the sole relation between the $\xi^i$, we can consider $\F_{p^{r-1}}$ as
the $\F_p$-module $V$, as above, with the generators $1,\xi,\ldots,\xi^{r-1}$,
and an expression \eqref{EqRep} of an element $a$ as sum of powers with as few
terms as possible corresponds to an admissible coordinate vector for $a$ as an
element of $V$.

Thus, as $\gcd(p,r)=1$, the result follows by Theorem~\ref{ThmG}. 
\proofend

\paragraph*{Proof of Theorem \ref{Thm2}.}
The nonzero $(p^{r-1}-1)/(2r)$th powers in $\F_{p^{r-1}}$ are exactly the 
$(2r)$th roots of unity in $\F_{p^{r-1}}$, and again $\F_{p^{r-1}}$ is
generated by a primitive $r$th root of unity. We consider the same module
$V$ as in the proof of Theorem \ref{Thm1}. Now, a representation of the form
\eqref{EqRep} with a minimal number of terms corresponds to an expression
$$
  a = \sum_{i=0}^{r-1} \pm v_i \xi^i,
$$
with $\sum |v_i|$ minimal; but this is the same as having
$$
  \|(\pm v_0, \ldots, \pm v_{r-1})\|_2
$$
minimal, where by the linear dependence of the $\xi^i$ we may add ${\bf
e}=(1,1,\ldots,1)$ if that reduces the norm. The problem is thus to
characterise admissible vectors for the norm $\|\cdot\|_2$. But this is done
in Theorem~\ref{ThmMR}.
\proofend

~

\noindent
Johann Radon Institute for \hfill Institut f\"ur Mathematik B\\
Computational and Applied Mathematics \hfill Technische Universit\"at Graz\\
Austrian Academy of Sciences \hfill Steyrergasse 30\\
Altenbergerstra\ss e 69 \hfill 8010 Graz, Austria\\
4040 Linz, Austria \hfill \\
arne.winterhof@oeaw.ac.at \hfill c.vandewoestijne@tugraz.at


\begin{thebibliography}{00}

\bibitem{brlipl} R.\,A. Brualdi, S. Litsyn, and V.\,S. Pless, Covering radius,
        in: Handbook of coding theory, 755--826, North-Holland, Amsterdam, 
        1998.

\bibitem{coholilo} G.\,D. Cohen, I.\,S. Honkala, S. Litsyn, and A. Lobstein,
        Covering Codes, Elsevier, Amsterdam, 1997.

\bibitem{CvdWImpl} C.E. van de Woestijne, Implementation of the results of
        the present paper in KASH 2.x, available for download from
	\url{http://www.opt.math.tugraz.at/~cvdwoest/leenorm.kash}.

\bibitem{gaso} C. Garcia and P. Sol\'{e},
        Diameter lower bounds for Waring graphs and multiloop networks,
        Discrete Math. 111 (1993), 257--261. 

\bibitem{he} T. Helleseth, On the covering radius of cyclic linear codes and 
        arithmetic codes, Discrete Appl. Math. 11 (1985), 157--173. 

\bibitem{Small} C. Small, Diagonal equations over large finite fields, Canad.
	J. Math. 36 (1984), 249--262.

\bibitem{Weil} A. Weil, Numbers of solutions of equations in finite
	fields, Bull. Amer. Math. Soc., 55 (1949), 497--508.

\bibitem{wi1} A. Winterhof, On Waring's problem in finite fields, 
        Acta Arith. 87 (1998), 171--177. 

\bibitem{wi2} A. Winterhof, A note on Waring's problem in finite fields,
        Acta Arith. 96 (2001), 365--368.

\end{thebibliography}

\end{document}